\documentclass[12pt]{amsart}
\usepackage{hyperref}
\usepackage[latin1]{inputenc}
\usepackage[all]{xy}

\newtheorem{theorem}{Theorem}[section]

\newtheorem{lemma}[theorem]{Lemma}
\newtheorem{corollary}[theorem]{Corollary}
\newtheorem{remark}[theorem]{Remark}
\newtheorem{definition}[theorem]{Definition}

\numberwithin{equation}{section}

\newcommand{\C}{{\mathbb C} }

\newcommand{\cE}{{\mathcal E} }
\newcommand{\cF}{{\mathcal F} }

\newcommand{\cK}{{\mathcal K} }

\newcommand{\cO}{{\mathcal O} }

\def\ol#1{{\overline{#1}}}
\def\ii{\sqrt{-1}}
\def\pt{\partial}

\def\ka{K\"ah\-ler}
\def\he{Her\-mite-Ein\-stein}

\begin{document}

\title{Vector bundles on Sasakian manifolds}

\author[I. Biswas]{Indranil Biswas}

\address{School of Mathematics, Tata Institute of Fundamental
Research, Homi Bhabha Road, Bombay 400005, India}

\email{indranil@math.tifr.res.in}

\author[G. Schumacher]{Georg Schumacher}

\address{Fachbereich Mathematik und Informatik,
Philipps-Universität Marburg, Lahnberge, Hans-Meerwein-Strasse, D-35032
Marburg, Germany}

\email{schumac@mathematik.uni-marburg.de}



\date{}

\begin{abstract}

We investigate the analog of holomorphic vector bundles in the context of
Sasakian manifolds.

\end{abstract}

\maketitle

\section{Introduction}

Just as the contact manifolds serve as substitutes for symplectic manifolds in
odd dimensions, the Sasakian manifolds are the odd dimensional counterparts of
K\"ahler manifolds. Although they are around for quite some time, in recent
years they were extensively studied. This resurgence is
substantially influenced by the
recently found relevance of Sasakian manifolds in string theory. C. Boyer and
K. Galicki published a series of papers investigating various differential
geometric aspects of Sasakian manifolds (see \cite{BGsusy} and references
therein). Sasakian manifolds appeared in string theory through the work of J.
Maldacena \cite{Mal}; see the papers of J. Sparks and references
therein for more details.

Our aim here is to investigate the analog of holomorphic vector bundles in the
context of Sasakian manifolds, reflecting the rich theory of holomorphic vector
bundles on K\"ahler manifolds, more precisely, the theory of moduli for
stable
Sasakian holomorphic vector bundles. We show the existence of the moduli space
in the category of reduced complex spaces, and construct a K\"ahler structure
on it.

Necessary prerequisites include a Hodge theory
related to both the Sasakian and the holomorphic aspects. We arrive at
a notion of a \he\ connection. Every stable Sasakian holomorphic vector bundle
possesses a unique such connection. The variation of the \he\ structures in a
holomorphic family is closely related to infinitesimal deformations of the
Sasakian holomorphic structures. We introduce an intrinsic \ka\ metric on the
moduli space, and show that, up to a numerical constant, it is  equal to
the Chern form of a determinant line bundle equipped with a Quillen
metric.

\section{Sasakian manifolds}

Let $(X,g)$ be a connected oriented smooth Riemannian manifold of dimension
$2n+1$, where $n$ is a nonnegative integer. Denote by  $\nabla$ the Levi-Civita
connection. We recall the following definition.

\begin{definition}{\cite[Definition-Theorem 10]{BGsusy}}\label{de:sasaki}
The pair $(X\, ,g)$ is called a \textit{Sasakian manifold} if any of the
following three equivalent conditions hold:
\begin{enumerate}
\item[(i)] There is a Killing vector field $\xi$ on $X$ of unit length such
    that the section
\begin{equation}\label{Phi}
\Phi \, \in\, C^\infty(X,\, TX\otimes (TX)^*)
\end{equation}
defined by $\Phi (v) \,=\, -\nabla_v\xi$ satisfies the identity
\begin{equation}\label{id.}
(\nabla_v \Phi) (w)\, =\, g(v\, ,w)\xi- g(\xi\, ,w)v
\end{equation}
for all $v\, ,w\,\in\, T_xX$ and all $x\,\in\, X$.

\item[(ii)] There is a Killing vector field $\xi$ on $X$ of unit length
    such that the Riemann curvature tensor $R$ of $(X\, ,g)$ satisfies the
    identity
$$
R(v\, ,\xi)w \, =\, g(\xi\, ,w)v- g(v\, ,w)\xi
$$
for all $v$ and $w$ as above.

\item[(iii)] The metric cone $({\mathbb R}_+\times X\, , dr^2 \oplus r^2g)$
    is K\"ahler.
\end{enumerate}
\end{definition}

A Killing vector field $\xi$ on $X$ of unit length satisfies condition (i), if
and only if it satisfies condition (ii). Given a Killing vector field $\xi$ of
unit length satisfying condition (i), the Kähler structure on ${\mathbb
R}_+\times X$ asserted in statement (iii) is constructed as follows. Let $F$ be
the distribution of $X$ of rank $2n$ given by the orthogonal complement of
$\xi$. The homomorphism $\Phi$ (defined in \eqref{Phi}) preserves $F$, and
furthermore,
\begin{equation}\label{e2}
(\Phi\vert_F)^2 \, =\, -\text{Id}_F\, .
\end{equation}
Let $\widetilde{J}$ be the almost complex structure on ${\mathbb R}_+
\times X$ defined by the following conditions:
\begin{eqnarray}
\widetilde{J}\vert_F\,=\, \Phi\vert_F\\
\widetilde{J}\left(\frac{d}{dr}\right) \,=\, \xi\\
\widetilde{J}(\xi) \,=\, -\frac{d}{dr}.
\end{eqnarray}
This almost complex structure is integrable. The Riemannian metric
$dr^2 \oplus r^2g$ on ${\mathbb R}_+\times X$ is K\"ahler with respect
to the complex structure $\widetilde{J}$.

Conversely, if the metric cone $({\mathbb R}_+\times X\, , dr^2 \oplus r^2g)$
is K\"ahler, then consider the vector field on ${\mathbb R}_+\times X$ given by
$J(\frac{d}{dr})$, where $J$ is the almost complex
structure on ${\mathbb R}_+\times X$. The vector field $\xi$ on $X$ obtained by
restricting this vector field to $\{1\}\times X\, =\, X$
satisfies condition $(i)$, and hence condition $(ii)$, in
Definition~\ref{de:sasaki}.

We will consider the vector field $\xi$ (or equivalently, the K\"ahler
structure on ${\mathbb R}_+\times X$) as part of the definition of a Sasakian
manifold.

Let $X$ be a smooth oriented Riemannian
manifold of dimension $2n+1$ and $F\, \subset\, TX$
an oriented smooth distribution of rank $2n$. The quotient map
$$
TX \, \longrightarrow\, TX/F\, =: \, N
$$
defines a smooth one-form on $X$
\begin{equation}\label{e1}
\omega\, \in\, C^\infty(X, T^*X\otimes N)
\end{equation}
with values in the line bundle $N$. Since $X$ is oriented, the orientation of
$F$ induces an orientation of the normal bundle $N$. Therefore, $N$ has a
canonical smooth section given by the positively oriented vectors of unit
length in the fibers of $N$. Consequently, the form $\omega$ in \eqref{e1}
gives a nowhere vanishing smooth one-form on $X$. This one-form will also be
denoted by $\omega$. The distribution $F$ is said to be a \textit{contact
structure} on $X$ if the $(2n+1)$-form $(d\omega)^n\wedge \omega$ is nowhere
vanishing.

\begin{remark}\label{rem1}
{\rm The distribution $F$ is integrable, if it satisfies the Frobenius
condition which says that the one-forms $\omega$ satisfy the condition
$(d\omega)\wedge \omega\, =\, 0$. Therefore, a contact structure $F$ is not
integrable unless $n\, =\, 0$.}
\end{remark}

Now let $(X\, , g\, , \xi)$ be a Sasakian manifold. The distribution $F$ on $X$
of rank $2n$ given by the orthogonal complement of the Killing vector field
$\xi$ defines a contact structure on $X$. We note that the corresponding
one-form $\omega$ is the dual of $\xi$ with respect to the metric $g$. From the
condition that $(d\omega)^n\wedge \omega$ is nowhere vanishing it follows that
the restriction of $d\omega$ to $F$ is fiberwise nondegenerate.

\begin{lemma}\label{lem0}
For all $x\,\in\, X$ and all $v\, ,w\, \in\, F_x$,
\begin{equation}\label{eq:deomega}
d\omega (v\, ,w) \, =\, - g(\Phi (v)\, ,w)\, ,
\end{equation}
where $\Phi$ is defined in \eqref{Phi}.
\end{lemma}

\begin{proof}
{}From the definition of $\Phi$,
$$
- g(\Phi (v)\, ,w)\, =\, g(\nabla_v\xi\, ,w)\, .
$$
Since $\xi$ is a Killing vector field,
\begin{equation}\label{KVF}
g(\nabla_v\xi\, ,w) + g(\nabla_w\xi\, ,v)\, =\, 0\, .
\end{equation}

Extend $v$ and $w$ to smooth sections $\widetilde{v}$ and
$\widetilde{w}$ of $F$. Since $\widetilde{w}$ is orthogonal to $\xi$,
$$
g(\nabla_v\xi\, ,w)\, =\, - g(\xi\, ,\nabla_v \widetilde{w})\, .
$$
Using \eqref{KVF},
$$
g(\nabla_v\xi\, ,w)\,=\, - g(\nabla_w\xi\, ,v)\,=\,
g(\xi\, ,\nabla_w \widetilde{v})
$$
because $\widetilde{v}$ is also orthogonal to $\xi$. Therefore,
$$
- g(\Phi (v)\, ,w)\, =\, \frac{1}{2}(- g(\xi\, ,\nabla_v
\widetilde{w})+
g(\xi\, ,\nabla_w \widetilde{v}))\, =\, -\frac{1}{2}g(\xi\, ,
[\widetilde{v}\, ,\widetilde{w}])\, .
$$
But $-\frac{1}{2}g(\xi\, ,
[\widetilde{v}\, ,\widetilde{w}])\, =\, d\omega(v\, ,w)$
because both $\widetilde{v}$ and $\widetilde{w}$ are orthogonal
to $\xi$.
\end{proof}

The line subbundle of $TX$ generated by $\xi$ will
be denoted by $N$.

A compact connected Sasakian manifold $(X\, ,g\, , \xi)$ is called
\textit{quasi-regular} if all the orbits of the unit vector field $\xi$ are
closed.

\section{Holomorphic vector bundles}

\subsection{Differential forms and Hodge type
decomposition}\label{se3.1}

Let $(X\, ,g\, , \xi)$ be a Sasakian manifold. Let $E$ be a $C^\infty$ complex
vector bundle over $X$.

\begin{definition}\label{de:transv}
A locally defined complex differential form $\alpha$ of class $C^\infty$ on $X$
will be called \textit{transversal} if
$$
i_\xi \alpha\, =\, 0 \,\,\,~~\text{~~and~~}\,\,\,~~
i_\xi d\alpha\, =\, 0\, ,
$$
where $i_\xi$ is the contraction with the Killing vector field $\xi$.
\end{definition}

Since $d^2\, =\,0$, it follows immediately that for any transversal
differential form $\alpha$ also $d\alpha$ is transversal.

Let $L_\xi$ denote the Lie derivative with respect to $\xi$. Since
$$
L_\xi \alpha\, =\, di_\xi \alpha + i_\xi d\alpha\, ,
$$
we conclude that $\alpha$ is transversal, if and only if
$$
i_\xi \alpha\, =\, 0\, =\, L_\xi \alpha\, .
$$
Let $\alpha$ be a locally defined transversal complex $d$-form. Take any point
$x\, \in\, X$. Since $i_\xi \alpha\, =\, 0$, the evaluation of $\alpha(x)$ on
$\bigwedge^d T_xX \otimes_{\mathbb R}{\mathbb C}$ is determined by the
evaluation of $\alpha(x)$ on the subspace
$$
\bigwedge\nolimits^d F_x\otimes_{\mathbb R}{\mathbb C}\, \subset
\, \bigwedge\nolimits^d T_xX\otimes_{\mathbb R}{\mathbb C}\, ,
$$
where $F_x$ as before is the orthogonal complement of the line $\xi(x)\,
\subset\, T_xX$.

We recall that
\begin{equation}\label{e3}
(\Phi\vert_{F_x})^2 \, =\, -\text{Id}_{F_x}
\end{equation}
(see \eqref{e2}). Let
$$
\Phi^{\mathbb C}_x\, :=\, \Phi\vert_{F_x}\otimes_{\mathbb R}
{\mathbb C}\, :\, F_x\otimes_{\mathbb R}{\mathbb C}
\, \longrightarrow\, F_x\otimes_{\mathbb R}{\mathbb C}
$$
be the complexification of the
automorphism $\Phi\vert_{F_x}$. From \eqref{e3}
it follows that the eigenvalues of $\Phi^{\mathbb C}_x$ are $\pm\sqrt{-1}$.
We have a decomposition
$$
F_x\otimes_{\mathbb R}{\mathbb C}\, =\, F^{1,0}_x \oplus
F^{0,1}_x\, ,
$$
where $F^{1,0}_x$ and $F^{0,1}_x$ are the eigenspaces for the eigenvalues
$\sqrt{-1}$ and $-\sqrt{-1}$ respectively.

Now, for $p,q \, \geq\, 0$, define
$$
F^{p,q}_x\, :=\, (\bigwedge\nolimits^p  F^{1,0}_x)\otimes
(\bigwedge\nolimits^q  F^{0,1}_x)\, \subset\,
\bigwedge\nolimits^{p+q}  F_x\otimes_{\mathbb R}{\mathbb C}\, .
$$
Therefore, we have a decomposition
\begin{equation}\label{e4}
\bigwedge\nolimits^{d}  F_x\otimes_{\mathbb R}{\mathbb C}
\, =\, \bigoplus_{i=0}^d F^{i,d-i}_x
\end{equation}
for all $d\, \geq\, 0$.

Take a transversal complex $d$-form $\alpha$ defined over an open subset $U\,
\subset\,X$ such that $\alpha$ is nonzero at some point. We will say that
$\alpha$ is of {\it type} $(a,d-a)$, if the evaluation of $\alpha$ on
$F^{i,d-i}_x$ vanishes for all $i\, \not=\, a$ and all $x\, \in\, U$.

Let $F^{i,d-i}$ denote the $C^\infty$ subbundle of the vector bundle
$\bigwedge\nolimits^{d} F\otimes_{\mathbb R}{\mathbb C}$ defined by
the
condition that the fiber of $F^{i,d-i}$ over each point $x\, \in\, X$ is
$F^{i,d-i}_x$. The pointwise decomposition in \eqref{e4} gives the following
$C^\infty$ decomposition into a direct sum of vector bundles
\begin{equation}\label{e4p}
\bigwedge\nolimits^{d}  F\otimes_{\mathbb R}{\mathbb C}
\, =\, \bigoplus_{i=0}^d F^{i,d-i}\, .
\end{equation}

Now the decomposition $TX\, =\, F\oplus ({\mathbb R}\cdot\xi)$ induces a
decomposition
$$
\bigwedge\nolimits^d TX\otimes_{\mathbb R}{\mathbb C}
\,=\,
\bigwedge\nolimits^d F_x\otimes_{\mathbb R}{\mathbb C}
\oplus (\xi\otimes \bigwedge\nolimits^{d-1}
F_x\otimes_{\mathbb R}{\mathbb  C})
$$
for each $d\, \geq\, 0$. Combining this decomposition with the
decomposition in \eqref{e4p}, we have
\begin{equation}\label{p.f.}
\bigwedge\nolimits^d TX\otimes_{\mathbb R}{\mathbb C}
\,= \,\bigoplus_{i=0}^d F^{i,d-i}\oplus
\left(\xi\otimes \left(\bigoplus_{j=0}^{d-1} F^{j,d-j-1}\right)\right)\, .
\end{equation}

Consider the projection
$$
\bigwedge\nolimits^d TX\otimes_{\mathbb R}{\mathbb C}
\, \longrightarrow\, F^{i,d-i}
$$
obtained from the decomposition in \eqref{p.f.}. Its dual is an injective
homomorphism
\begin{equation}\label{p.f2.}
f_{d,i}\, :\, (F^{i,d-i})^*\, \longrightarrow\,
\bigwedge\nolimits^d T^*X\otimes_{\mathbb R}{\mathbb C}
\end{equation}
for each $d\, \geq\, 0$ and $i\leq d$. A nonzero transversal complex $d$-form
$\alpha$ is of type $(a,d-a)$, if and only if $\alpha$ is of the form
$f_{d,a}(\alpha')$ for some section $\alpha'$ of $(F^{a,d-a})^*$, where
$f_{d,a}$ is the injective homomorphism in \eqref{p.f2.}.

For any $i$ with $0\leq i\leq d$, consider the inclusion of $F^{i,d-i} \,
\hookrightarrow\, \bigwedge\nolimits^d TX\otimes_{\mathbb R} {\mathbb C}$ in
\eqref{p.f.}. Its dual is a projection
\begin{equation}\label{p.p.}
\phi_{d,i}\, :\, \bigwedge\nolimits^d T^*X\otimes_{\mathbb R}
{\mathbb C}\, \longrightarrow\, (F^{i,d-i})^*\, .
\end{equation}

Let $\alpha$ be a transversal differential form on $X$ of degree $d$. Take any
$i$ with $0\leq i\leq d$. Let
\begin{equation}\label{p.p2.}
\alpha_i\, :=\,  f_{d,i}\circ \phi_{d,i}(\alpha)
\end{equation}
be the differential form on $X$, where $f_{d,i}$ and $\phi_{d,i}$ is
constructed in \eqref{p.f2.} and \eqref{p.p.} respectively.

\begin{lemma}\label{lem1}
The differential form $\alpha_i$ in \eqref{p.p2.} is transversal.
\end{lemma}

\begin{proof}
{}From the construction of $\alpha_i$ is follows immediately that $i_\xi
\alpha_i\, =\, 0$.

To prove that $L_\xi \alpha_i\, =\, 0$, consider the identity in
\eqref{id.}. In it,
set $v\, =\,\xi$, and set $w$ to be a smooth section of $F\, =\,
\xi^{\perp}$. The identity in \eqref{id.} immediately implies that
\begin{equation}\label{id2}
(\nabla_\xi \Phi) (w)\, =\, 0\, .
\end{equation}

The automorphism $\Phi$ acts on $T^*X\otimes_{\mathbb R}
{\mathbb C}$. From
\eqref{id2} it follows that the bundles consisting of eigenvectors associated
to the action of $\Phi$ are preserved by $\xi$. In particular,
$$
L_\xi (f_{1,0}\circ \phi_{1,0})\, =\, 0\,=\,
L_\xi (f_{1,1}\circ \phi_{1,1})
$$
because $f_{1,0}\circ \phi_{1,0}$ and $f_{1,1}\circ \phi_{1,1}$
are projections to the eigenbundles for the eigenvalues
$\sqrt{-1}$ and $-\sqrt{-1}$ respectively. Since
$f_{d,j}\circ \phi_{d,j}$ are constructed from $f_{1,0}\circ
\phi_{1,0}$ and $f_{1,1}\circ \phi_{1,1}$, we now have,
\begin{equation}\label{e2.}
L_\xi (f_{d,j}\circ \phi_{d,j})\, =\, 0
\end{equation}
for all $j$. In other words, the decomposition
$$
\text{Id} \,=\, \bigoplus_{j=1}^d f_{d,j}\circ \phi_{d,j}
$$
is preserved by the flow on $X$ defined by $\xi$.

{}From \eqref{e2.} we have
$$
L_\xi (f_{d,i}\circ \phi_{d,i}(\alpha)) \, =\,
f_{d,i}\circ \phi_{d,i} (L_\xi\alpha)\, .
$$
But $L_\xi\alpha\, =\, 0$, because $\alpha$ is transversal. This completes the
proof of the lemma.
\end{proof}

Using the trivialization of $N$ given by $\xi$, the $N$-valued one-form
$\omega$ becomes a real-valued one-form on $X$. Now we consider the nowhere
degenerate form $d\omega|F$. Using the orthogonal projection $TX\,
\longrightarrow\, F$, this defines a two-form on $X$ which is $d$-closed, but
not exact in general. Then form $d\omega|F$ is clearly transversal. Now Lemma
\ref{lem0} has the following corollary:

\begin{corollary}\label{kaehler}
The transversal form $d\omega|F$ is of type $(1\, ,1)$.
\end{corollary}

The formal adjoint of the multiplication operator
\begin{eqnarray*}
L\,:\,(F^{p,q})^*&\longrightarrow & (F^{p+1,q+1})^*\\
\eta &\longmapsto & \eta \wedge (d\omega|F)
\end{eqnarray*}
is denoted by
\begin{eqnarray}\label{co}
\Lambda_\omega \,:\, (F^{p+1,q+1})^*&\longrightarrow  &
(F^{p,q})^*\, .
\end{eqnarray}

\subsection{Partial connections}\label{sc:partconn} Let
\begin{equation}\label{S}
S\, \subset\, TX\otimes_{\mathbb R} {\mathbb C}
\end{equation}
be a subbundle of positive rank, which is integrable. In other words, smooth
sections of $S$ are closed under the operation of the Lie bracket. Let $E$ be a
$C^\infty$ complex vector bundle over $X$. Let
\begin{equation}\label{S1}
q_S\, :\, T^*X\otimes_{\mathbb R} {\mathbb C}\, =\, (TX
\otimes_{\mathbb R} {\mathbb C})^*\, \longrightarrow\, S^*
\end{equation}
be the dual of the inclusion map of $S$ in $TX\otimes_{\mathbb R} {\mathbb C}$.

A \textit{partial connection on $E$ in the direction of} $S$ is a $C^\infty$
differential operator
\begin{equation}\label{S3}
D\, :\, E \, \longrightarrow\, S^*\otimes E
\end{equation}
satisfying the Leibniz condition. The Leibniz condition says that for a smooth
section $s$ of $E$ and a smooth function $f$ on $X$,
$$
D(fs) \,=\, fD(s) + q_S(df)\otimes s\,
$$
holds, where $q_S$ is the projection in \eqref{S1}.

The condition that $D$ satisfies the Leibniz condition immediately implies that
the order of the differential operator $D$ is exactly one. We recall
that the
symbol of a $C^\infty$ differential operator
$$
D'\, :\, A\, \longrightarrow\, B
$$
of order one is a $C^\infty$ section of $A^*\otimes B\otimes
(TX\otimes_{\mathbb R} {\mathbb C})$.

Let $D$ be a partial connection on $E$ in the direction of $S$.
The symbol of the first order differential operator $D$ coincides
with $\text{Id}_E\otimes \widehat{q}_S$, where $\widehat{q}_S$ is
the smooth section of $S^*\otimes TX\otimes_{\mathbb R} {\mathbb
C}$ given by $q_S$ in \eqref{S1}. Indeed,
this follows immediately from the fact that $D$ satisfies
the Leibniz identity.

Since the distribution $S$ is integrable, smooth sections of $\ker(q_S)$ are
closed under the exterior derivation. Therefore, we have an induced exterior
derivation on the smooth sections of $S^*$
\begin{equation}\label{ed}
\widehat{d}\, :\, S^*\, \longrightarrow\, \bigwedge\nolimits^2
S ^ *
\end{equation}
which is a differential operator of order one.

Let $D$ be a partial connection on $E$ in the direction of $S$. Consider the
differential operator
$$
D_1\, :\, S^*\otimes E\, \longrightarrow\, (\bigwedge\nolimits^2 S^*)
\otimes E
$$
defined by
$$
D_1(\theta\otimes s) \, =\, \widehat{d}(\theta)\otimes s
-\theta\wedge D(s)\, ,
$$
where $\widehat{d}$ is constructed in \eqref{ed}. The composition
\begin{equation}\label{ed2}
E \, \stackrel{D}{\longrightarrow}\, S^*\otimes E
\, \stackrel{D_1}{\longrightarrow}\, (\bigwedge\nolimits^2 S^*)
\otimes E
\end{equation}
is $C^\infty(X)$-linear. Therefore, the composition in \eqref{ed2} defines a
$C^\infty$ section
\begin{equation}\label{ed3}
{\mathcal K}(D)\, =\,
C^\infty(X, (\bigwedge\nolimits^2 S^*) \otimes E\otimes E^*)\, =\,
C^\infty(X, (\bigwedge\nolimits^2 S^*) \otimes \text{End}(E))\, .
\end{equation}

The section ${\mathcal K}(D)$ in \eqref{ed3} is called the \textit{curvature}
of $D$. If
$$
{\mathcal K}(D)\, =\, 0\, ,
$$
then the partial connection $D$ is called \textit{flat}.

\subsection{Holomorphic hermitian vector bundles}\label{se:holvec}

Let $(X\, ,g\, , \xi)$ be a Sasakian manifold. Let
\begin{equation}\label{w01}
\widetilde{F}^{0,1}\, :=\, F^{0,1}\oplus (\xi\otimes_{\mathbb R}
\mathbb C)\, \subset\, TX\otimes_{\mathbb R} {\mathbb C}
\end{equation}
be the distribution, where $F^{0,1}$ is constructed in \eqref{e4p}.

\begin{lemma}\label{lem2}
The distribution $\widetilde{F}^{0,1}$ in \eqref{w01} is integrable.
\end{lemma}

\begin{proof}
Let $v$ and $w$ be smooth sections of $F$. Then from \eqref{id.},
$$
(\nabla_v \Phi) (w)\,=\, g(v\, ,w)\xi\, .
$$
In particular, this is a section of $\widetilde{F}^{0,1}$. Using this and
\eqref{id2} it follows that the distribution $\widetilde{F}^{0,1}$ is
integrable.
\end{proof}

\begin{definition}\label{def0}
A {\em Sasakian complex vector bundle} on the Sasakian manifold
$(X\, ,g\, ,
\xi)$ is a pair $(E\, ,D_0)$, where $E$ is a $C^\infty$ complex vector bundle
on $X$, and $D_0$ is a partial connection in the direction $\xi$.
\end{definition}

Since $N$ is a foliation on $X$ of dimension one, any partial
connection in the direction $N$ is automatically flat.

Note that the vector field $\xi$ is contained in the foliation distribution
$\widetilde{F}^{0,1}$ in \eqref{w01}. Therefore, any partial connection $D$ on
$E$ in the direction of $\widetilde{F}^{0,1}$ defines a partial connection on
$E$ in the direction of $N$.

\begin{definition}\label{def1}
 A  {\em holomorphic vector bundle} on the Sasakian manifold
    $(X\, ,g\, , \xi)$ is a pair $(\{E\, ,D_0\}\, ,D)$, where $(E\, , D_0)$
    is a $C^\infty$ Sasakian complex vector bundle on $X$, and $D$ is a
    flat partial connection on $E$ in the direction of
    $\widetilde{F}^{0,1}$ (constructed in \eqref{w01}) satisfying the
    compatibility condition that $D_0$ coincides with the partial
    connection on $E$, in the direction of $\xi$, defined by $D$.
\end{definition}

\begin{definition}\label{defp}
Let
     $(\{E\, ,D_{0}\}\, , D_{E})$ and $(\{E'\, ,D'_{0}\}\, ,
D_{E'})$ be two holomorphic vector bundles on $(X\, ,g\, , \xi)$. A fiberwise
$\mathbb C$-linear $C^\infty$ map
$$
\Psi \,: \,E'\,\longrightarrow\, E''
$$
is called {\em holomorphic}, if $\Psi$ intertwines $D_E$ and $D_{E'}$.
\end{definition}

\begin{remark}
{\rm Let $(\{E\, ,D_0\}\, ,D)$ be a Sasakian holomorphic vector bundle on $X$.
Then the dual $E^*$ also has a natural structure of a Sasakian holomorphic
vector bundle. If $(\{E'\, ,D'_0\}\, ,D')$ is another Sasakian holomorphic
vector bundle, then $E\otimes E'$ also has the structure of a Sasakian
holomorphic vector bundle. In particular, $\text{End}(E)\, :=\, E\otimes E^*$
carries a natural structure of a Sasakian holomorphic vector bundle.}
\end{remark}

\begin{definition}\label{def2}
A {\em hermitian structure} on a
Sasakian complex vector bundle
$(E\, ,D_0)$ is a $C^\infty$ hermitian structure
on the complex vector bundle $E$ preserved by the
partial connection $D_0$.
\end{definition}

So a hermitian structure on a Sasakian holomorphic vector bundle
induces a hermitian structure on the dual Sasakian holomorphic vector
bundle. Similarly, hermitian structures on two
Sasakian holomorphic vector bundles induce a hermitian structure
on their tensor product.

Let $(\{E\, ,D_0\}\, ,D)$ be a Sasakian holomorphic
vector bundle on $X$ equipped with a hermitian metric $h$.
Let $D'$ be a connection on $E$ preserving $h$ such that
the partial connection on $E$ in the direction of $\xi$
induced by $D'$ coincides with $D_0$.
Let $s$ and $t$ be two locally defined smooth sections of $E$
which are flat with respect to the partial connection $D_0$.
We have
\begin{equation}\label{idc}
\phi_{1,1}(d(h(s\, ,t)))\, =\, h(\phi_{1,1}(D'(s))\, ,t)
+h(s\, ,\phi_{1,1}(D'(t)))\, ,
\end{equation}
where $\phi_{i,j}$ is defined in \eqref{p.p.}. Since
$$
F\otimes_{\mathbb R}{\mathbb C}\, =\,
F^{0,1}\oplus F^{1,0}\, =\,F^{0,1}\oplus \overline{F^{0,1}}\, ,
$$
it follows from \eqref{idc} that there is a unique connection $\nabla$ on the
complex vector bundle $E$ satisfying the following two conditions:
\begin{itemize}
\item $\nabla$ preserves $h$, and

\item the partial connection on $E$ in the direction of
$\widetilde{F}^{0,1}$ induced by $\nabla$ coincides with $D$.
\end{itemize}

Let $\cK(E,h)\, :=\, {\mathcal K}(\nabla)$ be the curvature of the
above mentioned connection $\nabla$. We observe that ${\mathcal
K}(D)$ is a section of $(F^{1,1})^*\otimes E \otimes E^*$.

The Chern forms $c_j(E,h)$ of the Sasakian holomorphic hermitian
vector bundle are defined using $\cK(E,h)$ in the usual way.
More precisely,
$$
\det \left({\rm id}_E +\frac{\sqrt{-1}}{2\pi}\cK(E,h)\right)\, =\,
\sum_{i\geq 0} c_i(E,h)\, .
$$
We note that each $c_i(E,h)$ is a closed form of degree $i$,
and $c_i(E,h)$ is of type $(i\, ,i)$.

Henceforth, all Sasakian manifold considered will be assumed to
be compact.

\begin{definition}\label{de:HE}
A Sasakian holomorphic, hermitian vector bundle $(E,h)$ is called
{\em \he}, if
\begin{equation}\label{eq:HE}
\ii\Lambda_\omega \cK(E,h) = \lambda\cdot {\rm id}_E
\end{equation}
for some constant real number $\lambda$.
\end{definition}
The constant $\lambda$ is determined a priori. We define
$$
{\rm vol} (X)\,:=\,\int_X (d\omega)^n\wedge\omega\, .
$$
Then
$$
\frac{1}{2\pi} \lambda \cdot {\rm rk}E \cdot {\rm vol}(E) \,=\, \int_X
c_1(E,h)
(d\omega)^n\wedge\omega\,=:\, \deg_\omega(E,h)\, ,
$$
where $\dim X\, =\, 2n+1$,
is actually independent of the choice of $h$.

\subsection{Stable Sasakian holomorphic vector bundles}\label{stablesasa}

We are in a position to define the stability of a Sasakian
holomorphic vector
bundle with respect to the contact form $\omega$. Stability is defined in terms
of coherent subsheaves. We first define the structure sheaf
$\cO_X$ to be the sheaf of smooth functions $f$ with the property
that $df$ lies in the image of $f_{1,1}$, where $f_{1,1}$ is
constructed in \eqref{p.f2.}.

\begin{definition}\label{de:sheaf}
A {\em coherent sheaf} on a Sasakian manifold $X$ is a sheaf of
$\cO_X$-modules, which is locally the cokernel of a morphism of Sasakian vector
bundles of the form
$$
\cO^{\oplus a}_X \, \longrightarrow\, \cO^{\oplus b}_X
$$
(see Definition \ref{defp}).
\end{definition}

Given a coherent Sasakian sheaf $\cE$, there is a Sasakian holomorphic line
bundle $\det(\cE)$ associated to it. Given any local resolution
$$
0\, \longrightarrow\, V_n\, \longrightarrow\, V_{n-1}\,
\longrightarrow\,
\cdots \, \longrightarrow\, V_0 \, \longrightarrow\, \cE
\longrightarrow\, 0
$$
over some open subset $U$, where each $V_{i}$ is a Sasakian holomorphic vector
bundle for, and all the homomorphisms are morphisms of Sasakian vector bundles,
the restriction of $\det(\cE)$ to $U$ is identified with $\otimes_{i=0}^n
\det(V_i)^{(-1)^i}$. These locally defined line bundles patch together in a
natural way to define the line bundle $\det(\cE)$ over $X$.

The \textit{degree} $\deg \cE$ of a {\em torsionfree} coherent Sasakian sheaf
is defined to be
$$
\deg \cE\, :=\, \int_X c_1(\det(\cE))
(d\omega)^{n-1}\wedge\omega
$$
(recall that $X$ is compact).

\begin{definition}\label{de:stable}
Let $E\,=\, ( \{ E \, , D_0 \}\, ,D)$ be a Sasakian holomorphic vector
bundle, and $\cE$ the associated locally free $\cO_X$-module. The bundle $E$ is
called {\em stable} if for any Sasakian coherent subsheaf 
$\cE'\,\subset\, \cE$
such that the quotient $\cE/\cE''$ is torsionfree of positive rank, the
inequality
$$
\frac{\deg(\cE')}{{\rm rk}(\cE')} \,<\, \frac{\deg(\cE)}{{\rm rk}(\cE)}
$$
holds.

The bundle $E$ is called {\em semistable}, if
$$
\frac{\deg(\cE')}{{\rm rk}(\cE')} \,\leq\, \frac{\deg(\cE)}{{\rm
rk}(\cE)}\, .
$$
\end{definition}

\subsection{Sasakian \he\ metrics and applications}

Let $(X\, ,g\, , \xi)$ be a quasi-regular Sasakian
manifold. From now onwards
we will assume quasi-regularity.

Recall that $F\, =\, N^\perp\, \subset\, TX$ is the orthogonal
complement of the one dimensional foliation generated by the
Killing vector field $\xi$. We will call $F$ the \textit{horizontal
distribution}. In \eqref{e2} we noted that $\Phi$ preserves the
horizontal distribution.

The closedness of the $\xi$-orbits amounts to an action of the circle group
$$
{\mathbb S}^1\, =\, \{z\, \in\, {\mathbb C}\, \mid\, \,~\,|z|\, =\,
1\}
$$
on $X$.
However, the non-integrability of the contact structure defined by the
horizontal distribution $F$ means that we can just consider local slices
transversal to $\xi$, which defines a complex orbifold structure on the
quotient $X/{\mathbb S}^1$. The Sasakian holomorphic structure for complex
vector bundles discussed in Section~\ref{se:holvec} now can be studied locally.
This is actually sufficient to show the existence of \he\ metrics on stable
bundles.

\textbf{Existence of \he\ structures.} \textit{Let $(X\, , g\, , \xi)$ be a
compact quasi-regular Sasakian manifold, and let $( \{E\, , D_0 \} \, , D)$ be
a stable Sasakian holomorphic vector bundle on it. Then there exists a \he\
structure on $E$. Furthermore, the corresponding connection is uniquely
determined.}

We indicate, why the now standard methods of Donaldson \cite{do} and
Uhlenbeck-Yau \cite{uhy} are applicable.

Given an initial hermitian metric $h_0$ on $E$ in the sense of
Definition~\ref{def2}, we consider the heat equation for a $1$-parameter family
$\{h_t\}_{t\geq 0}$ of Sasakian hermitian metrics with curvature $\cK(h_t)$
\begin{equation}\label{eq:heat}
\frac{\pt h_t}{\pt t} h_t^{-1}\,=\, \Lambda_\omega \ii \cK(h_t) -
\lambda\cdot {\rm id_E}\, .
\end{equation}
This equation is not parabolic on $X$. We note that the two-form
$(d\omega)\vert_F$ as well as the Riemannian metric $g$ are ${\mathbb
S}^1$-invariant. Hence the contraction operation $\Lambda_\omega$ in \eqref{co}
is invariant under the action of ${\mathbb S}^1$. With the restriction on
hermitian metrics in Definition~\ref{def2}, the above heat equation is
${\mathbb S}^1$-invariant, and it exists as a parabolic equation on any slice.
The maximum principle is applicable because of the compactness of $X$. The
usual $C^0$-estimates for the pointwise norm of the curvature follow from
$$
\frac{d}{dt}|\Lambda_\cK(h_t)|^2\, =\, \square |\Lambda_\cK(h_t)|^2 -
|\nabla \Lambda_\cK(h_t)|^2,
$$
where $\square$ must be defined in the sense of the horizontal distribution
$F$. The existence of solutions for all $t\geq 0$ is provided. It follows that
$$
\frac{d}{dt}\int_X |\Lambda_\cK(h_t)|^2(t) \,=\, - \int_X |\nabla
\Lambda_\cK(h_t)|^2(t)\, ,
$$
and integration over $t$ yields the existence of a sequence $t_j \to \infty$
such that
$$
\int_X |\nabla \Lambda_\cK(h_t)|^2(t_j) \,\longrightarrow\, 0\, .
$$
The final step, following Uhlenbeck and Yau \cite{uhy}, of finding
suitable gauge
transformations such that the transformed hermitian metrics converge to a
Sasakian \he\ metric is based upon purely local arguments. In our case, it can
be carried through on the local slices transversal to $\xi$.

Now we need the notion of a holomorphic family of holomorphic Sasakian vector
bundles.

Let $(X\, ,g\, ,\xi)$ be a compact quasi-regular Sasakian manifold, and let $S$
be a parameter space, which will be a complex manifold (later we will also
consider reduced complex spaces). We equip $S$ with the flat \ka\ form. We
extend the Killing field $\xi$ to a vector field on $X\times S$ using the
decomposition of $T(X\times S)$, and denote it by $\xi_S$. Extend the
horizontal distribution $F$ to a distribution $F_S$ on $X\times S$ of
codimension one by taking the direct sum of $F$ with $TS$. Like in
Section~\ref{se3.1} we can now define $\Lambda^p F_S$ and $F^{p,q}_S$. As in
\eqref{w01} define
\begin{equation}\label{ffc}
{\widetilde F}^{0,1}_S\,:=\, F^{0,1}_S\oplus \xi_S
\otimes_{\mathbb R}{\mathbb C}\, .
\end{equation}

\begin{definition}\label{de:desfamvb}
Let $(X\, ,g\, ,\xi)$ be a compact quasi-regular Sasakian manifold.
\begin{enumerate}
\item[(i)] A {\em family $\{(\cE_s\, ,D_{0,s})\}_{s\in S}$  of Sasakian
    complex vector bundles parameterized by $S$} is given by a 
complex
    vector bundle
     $\cE$ of class $C^\infty$ over $X\times S$ together with a partial
     connection $D_{0,\cE}$ in the direction of $\xi$. For $s\,\in\, S$,
we define
     $(\cE_s\, ,D_{0,s}) \,:=\, (\cE\, ,D_{0,\cE})|_{X\times\{s\}}$.

\item[(ii)] A {\em holomorphic family of Sasakian vector bundles} 
is a
    family of Sasakian complex vector bundles
    $\{(\cE_s\, ,D_{0,s})\}_{s\in S}$ (as in (i)) together with
a flat partial connection $D$ on $\cE$ in the direction on
$\widetilde{F}_S^{0,1}$ (defined in \eqref{ffc}) satisfying the
following condition: for all $s\, \in\, S$, the partial connection
on $\cE_s$ in the direction of $\xi$ defined by $D$ coincides with
$D_{0,s}$. The restriction of $D$ to $\cE_s$ will be
denoted by $D_s$.
\end{enumerate}
\end{definition}

We observe that the complex orbifold structure on $X/{\mathbb S}^1$ guarantees
that for holomorphic families of Sasakian holomorphic vector bundles, the
Schlessinger conditions of deformation theory hold, in particular we have
universal deformations for simple bundles rather than semi-universal ones. A
Sasakian holomorphic vector bundle is called \textit{simple}, if all of its
automorphisms are constant scalar multiplications. (In the above setting, the
base was assumed to be smooth. The usual approach would be to consider
differentiable families of complex vector bundles first, and use the
integrability condition for the holomorphic structures, which defines the base
space $S$ as a subset of a smooth local ambient space, in the last step. In
this way all relevant differential operators exist on the ambient space of the
base \cite{f-s}.)

\begin{theorem}\label{thm1}
Let $(X, g, \xi)$ be a compact quasi-regular Sasakian manifold. Then there
exists the moduli space of stable, Sasakian holomorphic vector bundles on $X$
in the category of Hausdorff complex spaces.
\end{theorem}

\begin{proof}
Since $X$ is compact and quasi-regular, any section of ${\mathcal O}_X$
(defined earlier) is a constant function. Therefore, all the eigenfunctions of
an endomorphism of a holomorphic Sasakian vector bundle on $X$ are constant
functions. By this argument any stable Sasakian holomorphic vector bundle is
simple. So the {\em universal} deformations of stable Sasakian holomorphic
vector bundles exist. The local patching of the base spaces $S$ of universal
deformations yields a (reduced) complex space. Over these we are given
holomorphic families of Sasakian holomorphic vector bundles equipped with \he\
metrics, which can be chosen to depend in a $C^\infty$ way on the parameter.

The Hausdorff property follows from the existence and uniqueness of \he\
connections:

Let $S$ and $S'$ be reduced complex spaces, and let
$\{E_s,h_s\}_{s\in S}$ and $\{E'_s,h'_s\}_{s\in S'}$
be families of Sasakian holomorphic
\he\ bundles. The functor
$$
Isom_{S\times S'}(E\times S', S \times E')\,\longrightarrow\,
((Sets))\, ,
$$
which assigns to any space $R\,\longrightarrow\, S\times S'$ the set of
holomorphic isomorphisms of the bundles pulled back to
$R$ (and whose
morphisms are defined in the obvious way) is representable by a complex space
$$
\psi \,:\, I \,\longrightarrow\, S\times S'
$$
together with a universal object. The fibers of
$\psi$ are either empty or torsors for $\C^*$. Because of stability,
$\psi$ defines
a $\C^*$-bundle over its (closed analytic) image. Given the choice of relative
\he\ metrics for both families, the functor of isometries is represented
by
$$
\psi^{HE}\,:\,I^{HE} \,\longrightarrow\, S\times S'\, ,
$$
where $I\,\supset\, I^{HE}\,\longrightarrow\, S\times S'$ is a principal
$\text{U}(1)$-bundle. Let
$$
\Gamma\,=\,\{ (s,s')\,\in\, S\times S',\mid\, E_s \sim E'_{s'} \text{
isometrically equivalent}\}\, .
$$
For any $(s\, ,s')\,\in\, S\times S'$ the corresponding \he\ metrics are
gauge equivalent.
Now $$\Gamma\, =\, \psi^{HE}(I^{HE})\,=\,\psi(I)$$ is a closed analytic
subspace of $S\times S'$, which shows
the properness of the equivalence relation given by $\Gamma$.
\end{proof}

\begin{remark}
{\rm Since the group of holomorphic isometries of \he\ bundles is
connected, the moduli space is the union
of deformation spaces of stable Sasakian bundles.}
\end{remark}

\subsection{Sasakian moduli metric for vector bundles}

In this section we will introduce a \ka\ structure on the moduli space of
stable Sasakian holomorphic vector bundles, which is
constructed using a certain
determinant line bundle, equipped with a Quillen metric.

We first introduce the Sasakian moduli form. The construction will be
functorial with respect to base change.

Let $\{(\cE_s\, ,D_{0,s}\}\, ,D_s )\}_{s\in S}$ be a holomorphic family of
holomorphic Sasakian vector bundles on $(X\, ,g\, ,\xi)$, equipped with a
$C^\infty$ family of \he\ metrics $\{h_s\}$, amounting to a hermitian metric on
the complex vector bundle $(\cE\, ,D_{0\cE})$. Let $\cK(\cE,h)$ be the
curvature form of the unique hermitian connection on $\cE$ compatible with its
holomorphic structure. Therefore, $\cK(\cE,h)$ is a smooth section of
$(F^{1,1}_S)^*\otimes \text{End}(\cE)$.

\begin{definition}
The {\em Sasakian moduli form} for
vector bundles is defined as the following fiber integral
\begin{equation}
\omega^{SB}\,=\, \int_{(X\times S)/S} {\rm tr}\left( \ii\cK(\cE,h)\wedge
\ii\cK(\cE,h)\right) \frac{(d\omega)^{n}}{n!}\omega\, .
\end{equation}
\end{definition}

It follows from its construction that the Sasakian moduli form is a $d$-closed
real $(1,1)$-form on the parameter space $S$. Using the above arguments, one
can see that it possesses, in the case where $S$ is reduced but possibly
singular, a local $\pt\ol\pt$-potential of class $C^\infty$.

We know that the curvature form of a hermitian, Sasakian holomorphic vector
bundle induces a real $(1,1)$-form, with values in the endomorphism bundle, in
horizontal direction. The same holds for the curvature form of a hermitian
vector bundle over $X\times S$, which defines a family of such objects. For
this reason, we are in a position to carry over the methods from \cite{gafa}
literally to the orbifold case of \he\ orbifold vector bundles over the
orbifold $X/{\mathbb S}^1$ (cf. also \cite{bai}). We arrive at the following:

\begin{theorem}\label{th:wp}
The Sasakian moduli form for holomorphic vector bundles $\omega^{SB}$ is a \ka\
form on the moduli spaces of Sasakian holomorphic stable vector bundles. It
possesses local $\pt\ol\pt$-potentials of class $C^\infty$ (with respect to
ambient smooth spaces).
\end{theorem}

Let $M$ be a complex projective manifold and $L$ a very ample
line bundle over $M$. Fix a positive hermitian structure on $L$.
Let $\omega$ denote the K\"ahler form on $M$ given by its curvature.
Define
$$
X\, :=\, \{v\,\in\, L\, \mid \, \Vert v\Vert \,=\,1\}\, .
$$
The group ${\mathbb S}^1$ has the natural multiplication action on $X$. The
hermitian connection on $L$ gives a splitting of $TX$. We define a Riemannian
structure on $X$ by assigning the $\omega$ in the horizontal direction; the
metric in the vertical direction is given by the standard metric on ${\mathbb
S}^1$. It is easy to see that the Riemannian manifold $X$ equipped with the
vector field $\xi$ that is given by the action of ${\mathbb S}^1$, is a
Sasakian manifold. We note that the embedding
$$
M\, \hookrightarrow\, P(H^0(M,\, L)^*)
$$
given by the linear system lifts to a diffeomorphism
$$
X\, \longrightarrow\, H^0(M,\, L)^*
$$
which is an isometry in the horizontal direction.

\begin{definition}
A Sasakian manifold $(X\, ,g\, ,\xi)$ is called {\em of projective
type} if it possesses a holomorphic hermitian Sasakian line bundle
$(L^S\, ,h^S)$ whose Chern form is a positive hermitian
form on $F^{1,0}$.
\end{definition}

Now we construct the determinant line bundle $\lambda$ on the base of a
universal deformation. The construction is based upon the Riemann-Roch formula
for hermitian vector bundles by Bismut, Gillet and Soul\'e \cite{bgs}, which
also holds for orbifolds (cf.\ \cite{ma}).

Let $({\mathcal E}\, ,h)$ be a family of holomorphic hermitian vector
bundles of rank $r$ over $X$ parametrized by $S$. Let
$$
\pi\, : \, (X/{\mathbb S}^1) \times S \,\longrightarrow\,
X/{\mathbb S}^1
$$
be the projection.
The following formula holds for Chern characters for elements in the
Grothendieck's $K$-group:
\begin{gather}\label{eq:ch}
\chi\,:=\, (r^2 - ch(\cE\otimes \cE^*, h) ) \cdot ch\left(((L,h) -
(L,h)^{-1})^{\otimes(n-1)}\right)\\ \hspace{1cm}
=\, 2^{n-1}
c_1^2(\cE,h)c_1(L,h)^{n-1} + \text{ higher order terms}\, .\nonumber
\end{gather}

The $(1\, ,1)$-component of the push forward $\pi_*\chi$ in \eqref{eq:ch} is
the $(1\, ,1)$-component of the fiber integral, hence it is
$\omega^{SB}$
up to a numerical constant. In view of the Riemann-Roch formula, \cite{bgs},
the corresponding determinant line bundle is
\begin{equation}\label{lambda}
\lambda\,:=\, \det R^\bullet \pi_*\left((\cO_{X/{\mathbb
S}^1}^{r^2}-\cE\otimes \cE^{-1})
\otimes(\pi^*L- \pi^*L^{-1} )^{\otimes(n-1)} \right)\, .
\end{equation}
We indicate, how the above determinant line bundle can be computed for a
locally free sheaf $\cF$ say  and a proper holomorphic map $\pi\,:\,
Z\,\longrightarrow\, S$: Take the direct images $R^q\pi_*\cF$ first. These
locally possess free resolutions. Like in Section~\ref{stablesasa}, we define a
determinant line bundle associated to the resolution, which only depends upon
the sheaf $R^q\pi_*\cF$. Finally the tensor power of these invertible sheaves
with alternating exponent defines
$$
\lambda(\cF)\,:=\, \underline{\underline \det}(\cF)\,:=\,
R^\bullet\pi_*\cF\, .
$$
For a
short exact sequence
$$
0\,\longrightarrow\, \cF' \,\longrightarrow\,  \cF\,\longrightarrow\,
\cF''\,\longrightarrow\, 0
$$
we have that
$\underline{\underline \det}(\cF)\,=\,\underline{\underline
\det}(\cF')\otimes\underline{\underline \det}(\cF'')$. This fact allows for a
generalization of the definition of determinant line bundles. If $\cF$ denotes
a coherent sheaf on $Z$, then $\cF$ can be replaced by a locally free
resolution in the algebraic case (and more generally by the corresponding
simplicial objects in the non-algebraic case). Furthermore, the definition can
be extended to elements of the relative $K$-groups. If we let $\cF$ stand for
an element of the relative $K$-group arising from hermitian vector bundles,
like in the case under consideration, and $\pi : Z \to S$ is a smooth \ka\
morphism with relative \ka\ form induced by a closed form $\omega_Z$ on the
total space, the main theorem of \cite{bgs} states the existence of a Quillen
metric $h^Q$ on $\lambda(\cF)$ such that its first Chern form equals
\begin{equation}\label{eq:bgsrr}
c_1(\lambda(\cF),h^Q)= - \left(\int_{Z/S} {\rm td}(Z/S,\omega_Z)\cdot {\rm ch}(\cF,h)\right)_{(1,1)}
\end{equation}

In this way $\lambda$ of \eqref{lambda} is being defined and equipped with a
Quillen metric. It can be verified immediately that in our case the right hand
side of \eqref{eq:bgsrr} is equal to the fiber integral of $\chi$, which was
defined in \eqref{eq:ch}. Now we have that the Chern form is up to a numerical
constant equal to the Sasakian moduli metric form in the sense of orbifold
structures. Observe that these arguments also hold in the case of an orbifold
structure on a singular moduli space. The necessary techniques are in
\cite[Appendix]{f-s}.

\begin{theorem}
The Sasakian moduli metric form for vector bundles is up to a numerical factor
equal to the Chern form of the determinant line bundle.
$$
c_1(\lambda\, ,h^Q)\,\simeq\, \omega^{SB}.
$$
\end{theorem}


\end{document}